%Paper: funct-an/9308002
%From: GUIDO@MAT.UTOVRM.IT
%Date: Mon,  9 AUG 93 20:59 GMT

%%%%%% FORMAT %%%%%%%%%%%%
\magnification\magstep1
\hoffset=0.5truecm
\voffset=0.5truecm
\hsize=15.8truecm
\vsize=23.truecm
\baselineskip=14pt plus0.1pt minus0.1pt \parindent=19pt
\lineskip=4pt\lineskiplimit=0.1pt      \parskip=0.1pt plus1pt

\def\np{\par\noindent}

\catcode`@=11
\def\quad@rato#1#2{{\vcenter{\vbox{
        \hrule height#2pt
        \hbox{\vrule width#2pt height#1pt \kern#1pt \vrule width#2pt}
        \hrule height#2pt} }}}
\def\quadratello{\mathchoice
\quad@rato5{.5}\quad@rato5{.5}\quad@rato{3.5}{.35}\quad@rato{2.5}{.25} }
\catcode`@=12

\font\sectionfont=cmbx10 scaled\magstep1
\def\proof#1{\medskip\noindent{\bf Proof{#1}.}\quad}
\def\endproof{\hfill$\quadratello$\bigskip}
\def\section #1\par{\vskip0pt plus.3\vsize\penalty-75
    \vskip0pt plus -.3\vsize\bigskip\bigskip
    \noindent{\sectionfont #1}\nobreak\smallskip\noindent}
\def\claim#1#2\par{\vskip.1in\medbreak\noindent{\bf #1.} {\sl #2}\par
    \ifdim\lastskip<\medskipamount\removelastskip\penalty55\medskip\fi}
\def\rmclaim#1#2\par{\vskip.1in\medbreak\noindent{\bf #1.} {#2}\par
    \ifdim\lastskip<\medskipamount\removelastskip\penalty55\medskip\fi}
\def\remark#1#2\par{\vskip.1in\medbreak\noindent{\bf #1 Remark.} {\rm #2}\par
    \ifdim\lastskip<\medskipamount\removelastskip\penalty55\medskip\fi}
\def\itm#1{\item{$(#1)$}}

\def\ca{{\cal A}}
\def\A{{\cal A}}

\def\cb{{\cal B}}
\def\cc{{\cal C}}
\def\cd{{\cal D}}
\def\D{{\cal D}}
\def\ce{{\cal E}}
\def\E{{\cal E}}

\def\ch{{\cal H}}
\def\H{{\cal H}}

\def\ck{{\cal K}}

\def\cam{{\cal M}}

\def\cs{{\cal S}}

\def\bc{{\bf C}}
\def\Co{{\bf C}}

\def\Na{{\bf N}}
\def\br{{\bf R}}
\def\Re{{\bf R}}

\def\a{\alpha}
\def\b{\beta}
  
\def\d{\delta}
\def\e{\varepsilon}

\def\l{\lambda}  \def\La{\Lambda}
\def\m{\mu}
\def\n{\nu}

\def\X{\Xi}
\def\p{\pi}
\def\r{\rho}
\def\s{\sigma}
\def\t{\tau}
\def\c{\chi}
\def\f{\varphi}

\def\imply{\Rightarrow}
\def\coimply{\Leftarrow}

\def\L#1{L^{#1}({\cal A},\tau)}
\def\Lh#1{L^{#1}({\cal A},\tau)_h}
\def\Li#1#2{L^#1\cap L^#2}
\def\Lih#1#2{(L^#1\cap L^#2)_h}
\def\Lip{{Lip(\Re)}}
\def\norm#1#2{\| #1 \|_#2}
\def\Th{\hat T}
\def\tr#1{{}^t#1}
\def\parmat#1{\left(\matrix{#1}\right)}
\def\C#1{{{\cal C}_0^{#1}(\Re)}}
\def\exp#1{ {\rm e}^{#1} }
\def\ov{\overline}

\centerline{\bf NON SYMMETRIC DIRICHLET FORMS}   \par
\centerline{\bf ON SEMIFINITE VON NEUMANN ALGEBRAS}   \par
\bigskip
\bigskip
\centerline {Daniele GUIDO${}^1$, Tommaso ISOLA${}^1$, Sergio
SCARLATTI${}^2$}  \par  \bigskip
\centerline{${}^1$Dipartimento di Matematica. Universit\`a di Roma
``Tor Vergata''.}  \par
\centerline{Via della Ricerca Scientifica. I-00133, Roma, Italy.}  \par
\centerline{${}^2$Dipartimento di Matematica. Universit\`a di L'Aquila.}  \par
\centerline{Via Vetoio, Coppito. I-67100, L'Aquila, Italy.}  \par
\bigskip
\centerline{guido@mat.utovrm.it}        \par
\centerline{isola@mat.utovrm.it}        \par
\centerline{scarlatti@vxscaq.cineca.it} \par
\bigskip
\centerline{July 1993}
\bigskip
 \noindent
 {\bf Abstract.} The theory of non symmetric Dirichlet forms is generalized to
the non abelian setting, also establishing the natural correspondences among
Dirichlet forms, sub-Markovian semigroups and sub-Markovian resolvents within
this context. Examples of non symmetric Dirichlet forms given by derivations on
Hilbert algebras are studied.
 \bigskip

 \section Introduction.\par

The theory of non commutative Dirichlet forms, which originated from
the pioneering examples of L.~Gross [G] and the general analysis of
S.~Albeverio
and R.~H\o egh-Krohn [AH] (see also [AHO]), has nowadays drawn
a renewed interest between researchers ([DL1], [DL2], [DR], [D3], [Sa], [GL]
and [Ci]). There are different reasons which, in our opinion, explain
(and justify) the recent activity in this area. On the one side the
presence of a feed-back effect due to the increasing ability showed
by the commutative theory in handling successfully analytic and probabilistic
problems during the last fifteen years ([AR], [AMR], [MR], [D2], and ref.
therein). On the other side the great recent development of other new branches
of mathematics such as A.~Connes' non commutative geometry ([Co] and ref.
therein) and  quantum probability ([Pa], [AW] and ref. therein) with which the
theory of non commutative Dirichlet forms can naturally fit in. Let us remark
that up to now all works on non commutative Dirichlet forms treated the
generalization to a non abelian setting of the symmetric classical theory (see
[F]).
 \par
 In this paper, we develop the general theory of non symmetric Dirichlet
forms on a semifinite von~Neumann algebra $\A$. This means that we study
sesquilinear forms on the Hilbert space $L^2(\A,\t)$, requiring only  the so
called ``weak sector condition'' , which, at the form level,
roughly  means the antisymmetric part of the form must be controlled by the
symmetric one. In this sense our work can be seen as a non commutative
extension
of the theory of Dirichlet forms as it has been recently presented in [MR],
where this condition is assumed from the very beginning. It is worthwhile to
notice that these authors are able to produce a large amount of examples of
Dirichlet forms (see [MR] Chap. II). Among their  examples let us quote the
following simple one: consider the form
 $$
\E(u,v):=\sum_{i,j=1}^n\int a_{ij}(x){\partial{u}\over\partial{x_i}}
{\partial{v}\over\partial{x_j}}dx\eqno(0.1)
 $$
 where $u$ and $v$ are ${\cal C}^\infty$ functions with compact support in an
open set of $\Re^n$. If the functions $a_{ij}(x)$ are locally summable on $U$,
the symmetric part $[\tilde a_{ij}(x)]$ of the matrix-valued function
$[a_{ij}(x)]$ is uniformly bounded from below by a positive constant and
the entries of the antisymmetric part $[\check a_{ij}(x)]$ are $L^\infty$
functions, it can be proven that the form
$\{\E,D(\E)\}$ is closable and its closure  is a Dirichlet form.
 As it will be explained in Section~5, a natural generalization of the
preceding example is given by the form
 $$
\E(x,y):=\sum_{i,j=1}^n (d_jx,a_{ij}d_iy)\eqno(0.2)
 $$
where $a_{ij}:=\d_{ij}+c_{ij}$ and $[c_{ij}]$ is an antisymmetric matrix whose
entries are in the center of $\A$.
 As a consequence of the theory developed in this paper, we are able to
prove that, if $d_i$ are closable derivations and the intersection of their
domains is dense, then the form in (0.2) gives rise to a Dirichlet form.
 \par
 To better illustrate our results, we recall that throughout this paper $\A$
is a von Neumann algebra with a faithful, normal, semifinite trace  $\tau$.
Forms, semigroups, resolvents etc. are defined on the complex Hilbert space
$\L2$, even though many of their properties and relations require the
 real Hilbert space $\Lh2$ and its underlying order structure in an essential
way.
 \par
 In Section 1 we collect some preliminary material taken from [MR] on the
relationships between coercive closed forms on a Hilbert space and strongly
continuous contraction resolvents (resp. semigroups and their  generators)
satisfying sector condition. The last paragraphs of the section recall the
essentials of I.E.~Segal's theory of non commutative $L^p$ spaces on $\A$ (see
[N], [Se], [St]).
 \par
 In section 2 we establish the correspondence between Dirichlet forms,
sub-Mar\-ko\-vian semigroups and sub-Mar\-ko\-vian resolvents, thus
generalizing
the results of S.~Albeverio and R.~H\o egh-Krohn ([AH], see also [DL1]) to the
non symmetric case. This is made adapting the non-symmetric abelian definitions
and results in [MR] to the non commutative (semifinite) case.
 \par
  Section 3 is devoted to the extension of some properties of
sub-Markovian semigroups, already studied in [DL1], to the non symmetric
context. In particular we prove that sub-Markovian semigroups may be extended
to
$L^p$ spaces and study a class of sub-Markovian semigroups on $\L\infty$,
showing a correspondence between such semigroups and those on
$\L2$. Finally we study the consequences of complete positivity for
semigroups and Dirichlet forms such as the contraction property for semigroups
on $\L\infty$.
 \par
 In Section 4 we study derivations on Hilbert algebras and, based on previous
results in [Sa] and [DL1], we prove that, for a closed derivation on a
Hilbert algebra, the self-adjoint part of its domain is closed under Lipschitz
functional calculus and the whole domain is closed under the modulus operation.
 We also show that the corresponding norm inequalities (see (4.4) and (4.6))
hold, $i.e.$ such a derivation is a Dirichlet derivation in the sense of
E.B.~Davies and J.M.~Lindsay [DL1].
 Moreover, a non-abelian chain rule holds for the ${\cal C}^1$ functional
calculus of a self-adjoint operator. We notice that $\d$ need not be a
$^*$-derivation for the previous results to hold. Finally we show how
derivations which are not $^*$-invariant give rise naturally to (non symmetric)
Dirichlet forms.
 \par
 In section 5 we prove a theorem which gives rise to new examples  of non
commutative Dirichlet forms (and related semigroups). These examples are  of
the  previously mentioned type. They were already studied in [DL1] in
the symmetric case: this  simply corresponds to requiring the antisymmetric
part $[c_{ij}]$ in (0.2) to vanish.
 \par
 Lastly let us mention that  these results may be useful in the context of
open quantum systems and quantum statistical mechanics which, as it is known,
represent a natural physical arena where these mathematical theories have found
interesting applications (see e.g. references in [AH] and [DL1]).

\section Section 1. Preliminaries. \par

In this section we first collect definitions and facts about strongly
continuous semigroups and related objects, referring to [MR] for
proofs and further results, and then definitions and facts about $L^p$
spaces on $\{\ca,\t\}$, a von Neumann algebra with a faithful
semifinite normal trace, referring to classic works of [Se], [N] and
[St] for more detailed analysis and proofs, and to [T] for the general theory
of
von~Neumann algebras. \bigskip

It is well known that there is a bijective correspondence between
strongly continuous contraction resolvents $\{G_\a\}_{\a>0}$ on a
Banach space $X$, strongly continuous contraction semigroups
$\{T_t\}_{t>0}$ on $X$, and closed, densely defined linear operators
$\{L,\cd(L)\}$ on $X$, with the properties that $(0,\infty)\subset
\r(L)$, and $\| \a(\a-L)^{-1} \| \leq1,\ \forall \a>0$. \np
These objects are related by
$$\eqalignno{
G_\a&=(\a-L)^{-1},\ \a>0 &\cr
G_\a x&=\int_0^\infty \exp{-\a t}T_txdt,\ x\in X &\cr
Lx &= \lim_{t\downarrow0} {T_t x - x\over t},\quad x\in
\cd(L):=\{ x\in X: \lim_{t\downarrow0} {T_t x - x\over t}
\ {\rm exists } \} &\cr
T_tx&=\lim_{\a\to\infty} T_t^{(\a)}x := \lim_{\a\to\infty}
\exp{-\a t}\sum_{n=0}^\infty {(t\a)^n \over n!}(\a
G_\a)^nx, \ x\in X. & (1.1)
}$$

Recall now the theory of coercive closed forms. \par
Let $\ch$ be a complex Hilbert space, $\ck\subset\ch$ a real vector
subspace s.t. $\ck+i\ck=\ch$ and $(x,y)\in\br,\ \forall x,y\in\ck$,
and denote with $M_h:= M\cap\ck$, the real
part of $M\subset\ch$, and with $x^*:=y-iz$,
the adjoint of $x=y+iz,\ y,z\in\ck$. \np
Let $\ce:
\cd(\ce)\times\cd(\ce)\to \bc$, where $\cd(\ce)$ is a subspace of
$\ch$, be a real-positive, sesquilinear form on $\ch$, that
is, $\forall x,y,z\in\cd(\ce)$, $\a,\b\in\bc$, \par
$$\eqalign{
\ce(x,\a y+\b z)&=\a\ce(x,y)+\b\ce(x,z) \cr
\ce(y,x)&=\overline{\ce(x,y)} \cr
\ce(x^*,y^*)&=\overline{\ce(x,y)}, \cr
}$$
and
 $$
\ce(x,x)\geq0\ ,\qquad x\in \cd(\ce)_h,
 $$
and denote by $\tilde\ce$ the symmetric part of $\ce$,
$$
\tilde\ce(x,y):=
{1\over2}[\ce(x,y)+\overline{\ce(y,x)}],\ x,y\in\cd(\ce),
$$
and by $\check\ce$ the antisymmetric part of $\ce$,
$$
\check\ce(x,y):=
{1\over2}[\ce(x,y)-\overline{\ce(y,x)}],\ x,y\in\cd(\ce).
$$
Finally, denote by $\ce_\a$, $\a\geq0$, the form
$\ce_\a(x,y):=\ce(x,y)+\a(x,y),\ \forall x,y\in\cd(\ce)$. \par

\claim{1.1 Definition} $\{\ce,\cd(\ce)\}$ is said to satisfy the weak
sector condition if $\exists K>0$ s.t. $|\ce_1(x,y)| \leq
K\ce_1(x,x)^{1/2}\ce_1(y,y)^{1/2},\ x,y\in\cd(\ce)_h$.
\par

Notice that the above definition is equivalent to: $\exists  K'>0$ s.t.
$|\check\ce_1(x,y)| \leq K\ce_1(x,x)^{1/2}\ce_1(y,y)^{1/2},\ x,y\in\cd(\ce)_h$.
\par

\claim{1.2 Definition} $\{\ce,\cd(\ce)\}$ is said a
coercive closed form on $\ch$ if
\itm{i} $\cd(\ce)$ is dense in $\ch$
\itm{ii} $\{\tilde\ce,\cd(\ce)\}$ is closed [$i.e.\
\{ \cd(\ce),\tilde\ce_1 \}$ is a Hilbert space]
\itm{iii} $\{\ce,\cd(\ce)\}$ is real-positive and satisfies the
weak sector condition.
\par

\claim{1.3 Definition} A positive linear operator $\{A,\cd(A)\}$ on
$\ch$ is said to satisfy the sector condition if $\exists K>0$ s.t.
$|(x,Ay)|\leq K(x,Ax)^{1/2}(y,Ay)^{1/2},\ x,y\in\cd(A)_h$.
\par

\claim{1.4 Theorem} There is a bijective correspondence between
coercive closed forms $\{\ce,\cd(\ce)\}$ and strongly continuous
contraction resolvents $\{G_\a\}_{\a>0}$ s.t. $G_\a$ satisfies the
sector condition for some (hence for all) $\a>0$. \np
These objects are related by \np
$$
\ce_\a(x,G_\a y)=(x,y),\ x\in\cd(\ce),\ y\in\ch,
$$
and, if $L$ is the generator of $\{G_\a\}_{\a>0}$,
$$
\ce(x,y)=(x,-Ly),\ x\in\cd(\ce),\ y\in\cd(L),
$$
where $\cd(\ce)$ is the
completion of $\cd(L)$ w.r.t. $\tilde\ce_1^{1/2}$.
\par

\claim{1.5 Proposition} Let $\{\ce,\cd(\ce)\}$ be a coercive
closed form on $\ch$, and $\{G_\a\}_{\a>0}$, the associated
resolvent. Then, setting $\ce^{(\b)}(x,y):=\b(x,y-\b G_\b y),\
x,y\in\ch$, we get\np
\itm{i} $|\ce_1^{(\b)}(x,y)|\leq
(K+1)\ce_1(x,x)^{1/2}\ce_1^{(\b)}(y,y)^{1/2},\
x\in\cd(\ce),\ y\in\ch$
\itm{ii} Let $x\in\ch$. Then $x\in\cd(\ce)\iff \sup_{\b>0}
\ce^{(\b)}(x,x)<\infty$
\itm{iii} $\lim_{\b\to\infty} \ce^{(\b)}(x,y)=\ce(x,y),\
x,y\in\cd(\ce)$.
\par

\bigskip

Let $\A$ be a semifinite von Neumann algebra with a faithful
normal semifinite trace $\t$, and let $(\p,\H,\La)$ its GNS
representation. From now on $\A$ is identified with its
representation $\p(\A)$. \par
Let $\X$ be the collection of the closed, densely defined operators on $\H$
affiliated with $\A$. Then, by spectral theorem, $\forall x\in\X_h$,
$$
x=\int_{-\infty}^{+\infty}\l de_x(\l),
$$
where $e_x(E)\in\A$ for any Borel subset $E$ of $\br$, therefore
 $$
\n_x(E):=\t(e_x(E))
 $$
is a Borel measure on $\br$ and
 $$
\t(x):=\int_0^{+\infty}\l d\nu_x(\l),\qquad x\in\X_+
 $$
 is a faithful extension of $\t$ to $\X_+$.
\par
 Now, let us define, $\forall\ x\in\X$,
and for $p\in[1,\infty)$, $\|x\|_p:=\t(|x|^p)$.\np
Let, for each $p\in[1,+\infty)$,
 $$
\L p :=\{ x\in\X:\|x\|_p<+\infty \}
 $$
 and,
 $$
\tilde\A:=\{ x\in\X:\nu_{|x|}((\l,+\infty))<+\infty {\rm\ for\ some\ }
\l>0  \}.
 $$
Finally set $\L\infty:=\A$. It turns out that $\tilde\A$, equipped with strong
sense operations [Se] and with the topology of convergence in measure
([St], [N]), becomes a topological $*$-algebra, called the algebra of
$\t$-measurable operators. Moreover $\{\L p, \|\cdot\|_p\}$ is
a Banach subspace of $\tilde\A$, is linearly spanned by its positive
elements, and its norm satisfies $\|x\|_p=\|x^*\|_p$, for all
$x\in\L p$. Finally, observe that $\L2$ is a Hilbert space, with the
scalar product given by $(x,y):=\t(x^*y),\ x,y\in\L2$.\np

The basic properties of the $L^p$ spaces are:

\claim{1.6 Proposition} $(i)$ The trace $\t$ is extended to $\L1$ by
linearity, and
 $$
\eqalign{|\t(x)|&\leq\t(|x|)=\|x\|_1\cr
         \t(x^*)&=\overline{\t(x)}\cr}  \qquad x\in\L1
 $$
\item{$(ii)$} For each $x\in\X_h$ and for each Borel measurable function
$\f:\Re\to\Co$ one has
 $$
\|\f(x)\|_p=\|\f\|_{L^p(\Re,\n_x)}.
$$
In particular, if $\f\geq0$ or
$\f\in L^1(\Re,\n_x)$ we get
 $$
\t(\f(x))=\int \f d\n_x
 $$
\item{$(iii)$} $x,y\in\L p_+$, $x\leq y\imply\|x\|_p\leq\|y\|_p$.
\item{$(iv)$} $x,y\in\tilde\A$, $p,q,r\in[1,+\infty]$,
${1\over r}={1\over p}+{1\over q}\imply\|xy\|_r\leq\|x\|_p\|y\|_q$
\item{$(v)$}  Let $p\in[1,+\infty)$. Then for each $\psi\in\L p^*$ there exists
a unique element $x_\psi\in\L {p'}$, where ${1\over p}+{1\over p'}=1$,
such that  $\langle \psi, y\rangle=\t(x_\psi y)$ and
 %%%%%  $\|x_\psi\|_q=\|\psi\|_{\L p^*}$
$\psi\in\L p^*\to x_\psi\in\L {p'}$ is a Banach space isomorphism.
\item{$(vi)$} Let $p\in[1,+\infty]$, $x\in\L p$, $y\in\L {p'}$; then
 $$
\eqalign{\t(xy)&=\t(yx)\cr
\t(xy)&\geq0,\qquad x,y\geq0}
 $$
\item{$(vii)$} The weak$^*$ topology on $\L\infty$ given by the duality
$\langle\L1,\L\infty\rangle$ coincides with the $\s$-weak topology.
\par

 Let us denote with $Proj(\A)$ the set of self-adjoint idempotents
of $\A$ and with $\cs$ the following set
 $$
\cs:=\{x=\sum_{i=1}^n\l_ie_i\in\A:
n\in\Na,\ \l_i\in\Co,\ e_i\in Proj(\A)\cap\L1\}.
 $$

\claim{1.7 Proposition}
\itm{i} $\cs_+$ is dense in $L^p_+$ for $p\in[1,\infty)$ and weak$^*$
dense in $L^\infty_+$
\itm{ii} $x\in L^p,\ p\in[1,\infty]$ $\imply$ $xe_{|x|}(({1\over n},n))\to x$
in $L^p$
\itm{iii} Let $x\in\tilde\ca$; if $ax\in L^1$ and $\t(ax)=0,\ \forall
a\in\L1\cap\L\infty$ then $x=0$.
\par

\claim{1.8 Proposition} (Riesz-Thorin-Kunze interpolation) \np
Let $T:\cs\to\tilde\ca$ be a linear map satisfying
$$
\norm{Tx}{{q_i}}\leq
M_i\norm{x}{{p_i}},\ \forall x\in\cs,\ i=1,2,
$$
with $p_i,q_i\in[1,\infty]$ and $M_i>0$. If ${1\over p}:={t\over
p_1}+{1-t\over p_2}$, ${1\over q}:={t\over q_1}+{1-t\over q_2}$, where
$t\in(0,1)$, then, $\forall x\in\cs$, we get
$$
\norm{Tx}{q}\leq
M_1^{1-t}M_2^t\norm{x}{p}.
$$
\par

 \section Section 2. Markov semigroups and Dirichlet forms.
 \par

 In this section we give the basic definitions and prove the main theorems
which constitute the basis of the theory of non symmetric Dirichlet forms in a
non commutative setting. In our exposition we generalize to the non abelian
case results and techniques of Chap. I, Sec. 4 in [MR]. In particular, the
classical space of square integrable functions on a measure space is
replaced by the space of the operators affiliated to a von Neumann algebra $\A$
which are square integrable w.r.t. a normal, semifinite, faithful trace $\t$.
 \par
 \claim{2.1 Definition} $(i)$ A bounded linear operator $G$ on $\L2$ is called
{\it sub-Markovian} if
 $$
0\leq x\leq 1\imply 0\leq Gx\leq1\ ,\qquad\forall x\in\L2.
 $$
 A strongly continuous contraction resolvent $\{G_\a\}_{\a>0}$, resp. semigroup
$\{T_t\}_{t>0}$, is called sub-Markovian if all $\a G_\a$, $\a>0$, resp.
$T_t$,  $t>0$, are sub-Markovian.
 \item{$(ii)$} A closed densely defined operator $\{L,\cd(L)\}$ on $\L2$ is
called {\it Dirichlet operator} if $(Lx,(x-1)^+)\leq0$ for each $x\in{\cal
D}(L)_h$.
 \item{$(iii)$} A coercive closed form on $\L2$ is called a {\it
Dirichlet form} if, for all $x\in\D(\E)_h$, $x^+\wedge1\in\D(\E)$ and
 $$
\eqalign{
\E(x-x^+\wedge1,x+x^+\wedge1)&\geq0\cr
\E(x+x^+\wedge1,x-x^+\wedge1)&\geq0\cr}\eqno(2.1)
 $$
If only the first inequality in (2.1) holds, the form is called
$1/2$-Dirichlet.
 \par
 As in the classical case, if the form $\E$ is symmetric each of the two
inequalities in (2.1) is equivalent to the usual definition of Dirichlet form
(see e.g. [AH]).
 \par
 The following two theorems state the equivalence among the objects described
in Definition 2.1.
\claim{2.2 Theorem} Let $\{\E,\D(\E)\}$ be a coercive closed form on $\L2$
 with corresponding semigroup $\{T_t\}_{t>0}$, resolvent
$\{G_\a\}_{\a>0}$  and generator $\{L,\D(L)\}$. Then the following are
equivalent: \item{$(a)$} The form $\E$ is $1/2$-Dirichlet. \item{$(b)$} The
semigroup $\{T_t\}_{t>0}$ is sub-Markovian. \item{$(c)$} The resolvent
$\{G_\a\}_{\a>0}$ is sub-Markovian. \item{$(d)$} The generator $\{L,\D(L)\}$
is a Dirichlet operator.
 \par
 \claim{2.3 Theorem} Under the same hypotheses of the preceding theorem, the
following are equivalent:
 \item{$(a)$} The form $\E$ is Dirichlet.
 \item{$(b)$} The semigroups $\{T_t\}_{t>0}$ and $\{T_t^*\}_{t>0}$ are
sub-Markovian.
 \item{$(c)$} The resolvents $\{G_\a\}_{\a>0}$ and $\{G_\a^*\}_{\a>0}$ are
sub-Markovian. \item{$(d)$} The generators $L$ and  $L^*$ are Dirichlet
operators.
 \par
 The proof of the preceding theorems follows directly from propositions 2.6 and
2.7.
 \par
 \claim{2.4 Lemma} A bounded linear operator $G$ on $\L2$ is sub-Markovian
$iff$
 $$
\left\{\eqalign{x\geq0&\imply Gx\geq0\cr
x\leq1&\imply Gx\leq1\cr}\right.\qquad\forall x\in\L2
 $$
 \par
 \proof{} Sufficiency is true by definition. Now let
$x\in\L2$, $x\geq0$, and define $x_n:= x\wedge n$. Clearly $x_n\to x$ in
$\L2$, and therefore $Gx_n\to Gx$ in $\L2$ by continuity. Moreover $0\leq
{x_n\over n}\leq 1$ which implies $0\leq G({x_n\over n})\leq 1$, and therefore
$Gx_n\geq0$. Then, since the positive part of $\L2$ is closed, we get
$Gx\geq0$.
Finally let $x\in\L2$, $x\leq1$. If $x=x^+-x^-$ is the decomposition of $x$
into positive and negative part, we have $0\leq x^+\leq1$, $0\leq x^-$ and
therefore $0\leq Gx^+\leq 1$ by the sub-Markov property and $0\leq Gx^-$ by the
first part of this theorem. Then the thesis follows by linearity.
\endproof
\claim{2.5 Lemma} Let $\{x_n\}$ be a sequence converging to $x$ in
$\L2$ for which $0\leq x_n\leq 1$, $\forall n\in\Na$. Then $0\leq x\leq1$.
\par
\proof{} The fact that $x\geq 0$ follows because $\L2_+$ is norm closed.
Moreover, since $x_n$ converges weakly in $L^2$ and is uniformly bounded in
$\L\infty$, $x_n$ converges to $x$ weak$^*$ in
$\L\infty$ and therefore $\|x\|_\infty\leq1$. This implies $x\leq1$ because
$x$ is positive.
\endproof
\claim{2.6 Proposition} Let $\{G_\a\}_{\a>0}$ be a strongly
continuous contraction resolvent on $L^2(\A,\t)$ with corresponding
generator $L$ and semigroup $\{T_t\}_{t>0}$. Then the following are
equivalent: \item{$(i)$} $\{G_\a\}_{\a>0}$ is sub-Markovian.
\item{$(ii)$} $\{T_t\}_{t>0}$ is sub-Markovian.
\item{$(iii)$} $L$ is a Dirichlet operator.
\par
\proof{} $(i)\Rightarrow(ii)$:
Let $x\in L^2(\A,\t)$ and $0\leq x\leq1$. Then, for all $\b>0$,
$x_\b:=\b G_\b x$ is in $\cd(L)$ and $0\leq x_\b\leq1$ since  the resolvent
$G_\b$  is sub-Markovian, therefore, by formula (1.1)
and lemma 2.5, $0\leq T_tx_\b\leq1$. Moreover
$x_\b\to x$ in $\L2$ when $\b\to\infty$, therefore, again by
lemma 2.5,  $0\leq T_tx\leq1$, $i.e.$ $T_t$ is sub-Markovian.
\np
$(ii)\Rightarrow(iii)$: Let $x\in\Lh2$. Then
 $$
((x-1)^+,T_t(x-1)^+)\leq((x-1)^+,(x-1)^+)=((x-1)^+,(x-1))
 $$
by the Schwartz inequality and the fact that $(x-1)^+$ and $(x-1)^-$ are
orthogonal in $\L2$. Moreover $T_t(x\wedge1)\leq1$ by lemma 2.4. Therefore,
since $x=(x-1)^++x\wedge1$, we have
 $$
\eqalign{((x-1)^+,T_tx)&=((x-1)^+,T_t(x-1)^+)+((x-1)^+,T_t(x\wedge1))\cr
&\leq((x-1)^+,(x-1))+\t((x-1)^+)\cr
&=((x-1)^+,x).\cr}
 $$
Therefore we get
$$
((x-1)^+,Lx)=\lim_{t\downarrow0}{1\over t}((x-1)^+,T_tx-x)\leq0,
\qquad\forall x\in \cd(L).
 $$
$(iii)\imply(i)$. Let $x\in \Lh2$, $\a>0$ and $y:=\a G_\a
x$. We want to prove that if $0\leq x\leq1$ then $0\leq y\leq1$.
Indeed, for $x\leq1$, we have
 $$
\eqalign{\a((y-1)^+,y)&=((y-1)^+,\a y-Ly)+((y-1)^+,Ly)\cr
&\leq\a((y-1)^+,x)\leq\a\t((y-1)^+).\cr}
 $$
As a consequence,
 $$
\norm{(y-1)^+}{2}=((y-1)^+,y)-\t((y-1)^+)\leq0,
 $$
hence $y\leq1$. On the other hand, if $x\geq0$, then $-nx\leq1$ $\forall
n\in\Na$, therefore, by the previous result, $-ny\leq1$, $\forall n\in\Na$,
$i.e.$ $y\geq0$.
 \endproof
 \claim{2.7 Proposition} Let $\{\E,\D(\E)\}$ be a coercive closed form on $\L2$
with resolvent $\{G_\a\}_{\a>0}$. Then the following are equivalent:
\item{$(i)$} For all $x\in\D(\E)_h$ and $\a\geq0$, $x\wedge\a\in\D(\E)$ and
$\E(x-x\wedge\a,x\wedge\a)\geq0$.
\item{$(ii)$} For all $x\in\D(\E)_h$, $x^+\wedge1\in\D(\E)$ and
$\E(x-x^+\wedge1,x^+\wedge1)\geq0$.
\item{$(iii)$} $\E$ is a $1/2$-Dirichlet form.
\item{$(iv)$} $\{G_\a\}_{\a>0}$ is sub-Markovian.
\np
The analogous equivalences hold when $\{G_\a\}_{\a>0}$  is replaced by its
adjoint and $\E$ by the form $\E^\dagger(x,y):=\E(y,x)$.
\par
\proof{} $(i)\imply(ii)$. Let $x\in \D(\E)_h$, then, by $(i)$, we get
$x^-$, $x^+$, $x^+\wedge 1\in \D(\E)$. As a consequence
$$\eqalign{
\E(x-x^+\wedge1,x^+\wedge1)=&
\E(x^+-x^+\wedge1,x^+\wedge1)-\E(x^-,x^+\wedge1)\cr
\geq& -\E((x\wedge1)^-,(x\wedge1)^+).\cr}$$
Now for any $y\in \D(\E)_h$ we have, again by $(i)$,
$$
\E(y^-,y^+)=\E(y^+-y,y^+)=-\E((-y)-(-y)\wedge0,(-y)\wedge0)\leq0
$$
therefore $\E(x-x^+\wedge1,x^+\wedge1)\geq0$.\np
$(ii)\imply(iii)$. Since $\E$ is a real-positive sesquilinear form and
$(ii)$ holds, we get, for all $x\in \D(\E)_h$,
$$
\E(x-x^+\wedge1,x+x^+\wedge1)=\E(x-x^+\wedge1,x-x^+\wedge1)
+2\E(x-x^+\wedge1,x^+\wedge1)\geq0
$$
$(iii)\imply(iv)$. Let $y\in L^2(\A,\t)$, $0\leq y\leq 1$. We have
to show that $x:=\a G_\a y$ satisfies $0\leq x\leq 1$. Indeed
$$
\eqalign{
\|x&-x^+\wedge1\|^2+(x-x^+\wedge1,x^+\wedge1-y)=\cr
&=(x-x^+\wedge1,x-y)\cr
&=-{1\over\a}\E(x-x^+\wedge1,x)\cr
&=-{1\over2\a}\left(\E(x-x^+\wedge1,x+x^+\wedge1)
+\E(x-x^+\wedge1,x-x^+\wedge1)\right)\cr
&\leq0\cr
}\eqno(2.2)
$$
where the equality in the second line follows from theorem 1.4.
Let us introduce the functions
$f,g,h:\Re\to\Re$,
 $$\eqalign{
f(t)&=t\chi_{(-\infty,0]}(t),\cr
g(t)&=(t\wedge1)\chi_{[0,\infty)}(t),\cr
h(t)&=(t-1)\chi_{[1,\infty)}(t).\cr}
 $$
Then $fg\equiv0$, $gh\equiv h$, $g(x)=x^+\wedge1$ and
$x-g(x)=f(x)+h(x)$.
Therefore
 $$\eqalign{
(x-x^+\wedge1,x^+\wedge1-y)&=\t((x-x^+\wedge1)(x^+\wedge1-y))\cr
&=\t(f(x)(g(x)-y))+\t(h(x)(g(x)-y))\cr
&=\t((-f(x))y)+\t(h(x)(1-y))\geq0}\eqno(2.3)
$$
where we used proposition 1.6$(vi)$.
Finally equations (2.2) and (2.3) imply
$\|x-x^+\wedge1\|=0$, i.e. $0\leq x\leq1$.
\np
$(iv)\imply(i)$. Let $x\in \D(\E)_h$, $\a\geq0$. Now we prove that
$x\wedge\a\in\D(\E)$:
since $x=(x-\a)^++x\wedge\a$, it suffices to prove $(x-\a)^+\in\D(\E)$.
Recalling that, by proposition 1.5,
$\E^{(\b)}(y,z)=\b\ \t(y^*(z-\b G_\b z))$, for $y,z\in\L2$, we have
 $$
\eqalign{\E^{(\b)}((x-\a)^+,x\wedge\a)
&=\b\ \t((x-\a)^+(x\wedge\a))-\b\ \t((x-\a)^+\b G_\b(x\wedge\a))\cr
&\geq\a\b\ \t((x-\a)^+)-\a\b\ \t((x-\a)^+)=0\cr}\eqno(2.4)
 $$
where, since $x\wedge\a\leq\a$, the inequality in (2.4) follows from lemma
2.4, proposition 1.6$(vi)$ and the fact that $\b G_\b$ is sub-Markovian.
\np
Therefore,
 $$
\eqalign{ \E_1^{(\b)}((x-\a)^+&,(x-\a)^+)
=\E_1^{(\b)}((x-\a)^+,x-x\wedge\a)\cr
&=\E_1^{(\b)}((x-\a)^+,x)-\E_1^{(\b)}((x-\a)^+,x\wedge\a)\cr
&=\E_1^{(\b)}((x-\a)^+,x)-\E^{(\b)}((x-\a)^+,x\wedge\a)
                         -((x-\a)^+,x\wedge\a)\cr
&\leq\E_1^{(\b)}((x-\a)^+,x)\cr
&\leq(K+1)\E_1(x,x)^{1/2}\E_1^{(\b)}((x-\a)^+,(x-\a)^+)^{1/2}\cr}
 $$
where the last inequality follows from proposition 1.5$(i)$. As a consequence,
 $$
\E^{(\b)}((x-\a)^+,(x-\a)^+)\leq\E_1^{(\b)}((x-\a)^+,(x-\a)^+)
\leq(K+1)^2\E_1(x,x).\eqno(2.5)
 $$
Now proposition 1.5$(ii)$ and (2.5) imply $(x-\a)^+\in\D(\E)$.\np
Finally we prove that $\E(x-x\wedge\a,x\wedge\a)\geq0$: we have
 $$
\E^{(\b)}(x-x\wedge\a,x\wedge\a)=\E^{(\b)}((x-\a)^+,x\wedge\a)\geq0
 $$
by (2.4), hence the result follows by proposition 1.5$(iii)$.
\endproof
 We conclude this section with a theorem in which it is shown that a smooth
version of the definition of a Dirichlet form can be given. More precisely, the
{\it normal contraction} $x^+\wedge1$ in $(2.1)$ may be substituted by a
family of ${\cal C}^\infty$ contractions.
 \claim{ 2.8 Theorem} Let $\{\E,\D(\E)\}$ be a coercive closed form on $\L2$.
Then, the following are equivalent:
\item{$(i)$} $\E$ is a Dirichlet form.
\item{$(ii)$} For each $x\in\D(\E)_h$ there exists a family of functions
$\f_\e:\Re\to[-\e,1+\e]$, $\e\geq0$, such that
 \np
\hskip1.6cm{$(a)$} $\f_\e(t)=t$ for all $t\in[0,1]$.\np
\hskip1.6cm{$(b)$} $\f_\e$ is Lipschitz continuous with Lipschitz constant
1.\np
\hskip1.6cm{$(c)$} $\f_\e(x)\in\D(\E)$\np
\hskip1.6cm{$(d)$} $\liminf_{\e\to0}\E(x\mp\f_\e(x),x\pm\f_\e(x))\geq0$
 \par
\proof{} The implication $(i)\imply (ii)$ is trivial.
\np
$(ii)\imply (i)$.
Let us show that
 $$
\lim_{\e\to0}\|\f_\e(x)-x^+\wedge1\|_2 =0\ .\eqno(2.6)
 $$
Indeed, setting $\f_0(t):=(t\vee0)\wedge1$, we have
 $$
(\f_\e(x)-\f_0(x))^2=\int_{-\infty}^{+\infty}(\f_\e(\l)-\f_0(\l))^2de_\l
 $$
where $e_\l$ is the spectral family associated with $x$, therefore
 $$
\eqalign{\t((\f_\e(x)
&-\f_0(x))^2\cr
&=\t\left(\int_{-\infty}^{-\sqrt{\e}}(\f_\e(\l)  )^2de_\l+
          \int_{-\sqrt{\e}}^{0}      (\f_\e(\l)  )^2de_\l+
          \int_1^{+\infty}       (\f_\e(\l)-1)^2de_\l\right)\cr
&\leq\e^2\t\left(   \c_{(-\infty,-\sqrt{\e} ]}(x)\right)+
         \t\left(x^2\c_{[-\sqrt{\e} ,0      ]}(x)\right)+
     \e^2\t\left(   \c_{[      1,+\infty)}(x)\right).\cr}
 $$
 Then (2.6) follows by the following:
 $$
\t\left(\c_{[1,+\infty)}(x)\right)
\leq\t\left(x^2\c_{[1,+\infty)}(x)\right)
\leq\t(x^2)<\infty\ ,
 $$
 $$
\t\left(\c_{(-\infty,-\sqrt{\e}]}(x)\right)
\leq\t\left(\int_{-\infty}^{-\sqrt{\e}}{\l^2\over\e}de_\l\right)
\leq{1\over\e}\t(x^2)\ ,
 $$
 and
 $$
\t\left(x^2\c_{[-\sqrt{\e},0]}(x)\right)
=\t\left(\int_{-\sqrt{\e}}^0\l^2de_\l\right)
\to0,
 $$
when $\e\to0$, because $\m(E):=\t(\int_E\l^2de_\l)$ is a finite measure and
$\m(\{0\})=0$.
 \par
Summing the two inequalities in $(d)$ it follows
 $$
\limsup_{\e\to0}\E(\f_\e(x),\f_\e(x))\leq\E(x,x)
 $$
therefore, applying [MR, proposition 2.12] to  $\f_{\e}(x)$ and the form
$\E$, we get a sequence $\e_n\to0$ such that
$\f_{\e_n}(x)$ converges weakly in $\{\D(\E),{\tilde\E}_1\}$ to $x^+\wedge1$,
so that $x^+\wedge1\in\cd(\ce)$, and
 $$
\E(x^+\wedge1,x^+\wedge1)\leq\liminf_{n\to\infty}\E(\f_{\e_n}(x),\f_{\e_n}(x))
 $$
 Moreover, by the weak sector condition, the functional $\E_1(\cdot,x)$ is
continuous in $\{\D(\E)$, ${\tilde\E}_1\}$, therefore
 $$
\lim_{n\to\infty}\E(x,\f_{\e_n}(x))=\E(x,x^+\wedge1).
 $$
 Finally, we have
 $$
\eqalign{
\E(x\pm(x^+\wedge1),x\mp(x^+\wedge1))
&\geq\E(x,x)\mp\lim_{n\to\infty}\E(x,\f_{\e_n}(x))\cr
&\pm\lim_{n\to\infty}\E(\f_{\e_n}(x),x)
-\liminf_{n\to\infty}\E(\f_{\e_n}(x),\f_{\e_n}(x))=\cr
&=\limsup_{n\to\infty}\E(x\pm\f_{\e_n}(x),x\mp\f_{\e_n}(x))\geq0
}$$
where the last inequality follows by hypothesis $(d)$.
\endproof

\section Section  3. $L^p$ extensions of sub-Markovian
semigroups and complete positivity. \par

This section is devoted to the extension of some properties of
sub-Markovian semigroups, already studied in [DL], to the non symmetric
context. In particular we prove that sub-Markovian semigroups may be extended
to
$L^p$ spaces and study a class of sub-Markovian semigroups on $\L\infty$,
showing a correspondence between such semigroups and those on
$\L2$. Finally we exploit the consequences of complete positivity for
semigroups and Dirichlet forms such as the contraction property for semigroups
on $\L\infty$.
\par

\claim{3.1 Definition}
\itm{i} Let $M\in\cb(\L p)$, then we define $\tr{M}\in
\cb(\L {p'})$, where ${1\over p}+{1\over p'}=1$, as the unique linear
operator satisfying $(\tr{M}x,y)=(x,My),\ \forall x\in L^p,\ y\in L^{p'}$.
\np
\itm{ii} $M\in\cb(\L p)$, $p\in[1,\infty)$, is said a {\it sub-Markovian}
operator on $L^p$ if $x\in\L p$, $0\leq x\leq1$
$\imply$ $0\leq Mx\leq1$. \np
\itm{iii} $M\in\cb(\L \infty)$, is said a {\it sub-Markovian} operator on
$L^\infty$ if it is weak$^*$ continuous and $0\leq x\leq1$
$\imply$ $0\leq Mx\leq1$. \np
\itm{iv} $\{T_t\}_{t\geq0}\subset\cb(\L p)$, $p\in[1,\infty)$, is said a
{\it sub-Markovian semigroup} on $L^p$, if $T_t$ is a sub-Markovian
operator on $L^p$, for all $t>0$, and $T_t\to I,\ t\to0$, strongly on $\L
p$. \np \itm{v} $\{T_t\}_{t\geq0}\subset\cb(\L \infty)$ is said a {\it
sub-Markovian semigroup} on $L^\infty$, if $T_t$ is a sub-Markovian
operator on $L^\infty$, for all $t>0$, and $T_t\to I,\ t\to0$, weak$^*$ on
$\L \infty$.  \np
\par

\remark{3.2} Notice that the definition of sub-Markovian operator on
$L^2$ differs from that in the preceding sections in that we do not
require contractivity here.
\par

\claim {3.3 Proposition} Let $M$ and $M^*$ be sub-Markovian operators on
$L^2(\ca,\t)$. Then: \np
 \item{$(i)$} $\norm{Mx}{p}\leq 2\norm{x}{p}$, $\norm{M^*x}{p}\leq
2\norm{x}{p}$, $x\in L^p\cap L^2$, $p\in[1,\infty]$. \np
Let now $M^{(p)}$, resp. $M^{*(p)}$, be the unique continuous
extensions of $M|_{L^p\cap L^2}$, resp. $M^*|_{L^p\cap L^2}$, to
$L^p$, for $p\in[1,\infty)$. Then: \np
 \item{$(ii)$} $\tr{ (M^{(p)}) }=M^{*(p')}$, resp. $\tr{ (M^{*(p)})
}=M^{(p')}$,
for $p\in(1,\infty)$, and ${1\over p}+{1\over p'}=1$. \np
Set $M^{(\infty)} := \tr{ (M^{*(1)}) }$, resp. $M^{*(\infty)} :=
\tr{ (M^{(1)}) }$. Then: \np
 \item{$(iii)$} $M^{(p)}x=M^{(q)}x$, resp. $M^{*(p)}x=M^{*(q)}x$, for $x\in
L^p\cap L^q$, $p,q\in[1,\infty]$.
\par
\proof{}
$(i)$ Let $x\in\Lih{\infty}{2}$, then $\norm{x}{\infty}\leq1\iff
0\leq x_\pm\leq1\imply 0\leq Mx_\pm\leq1$ so that
$\norm{Mx}{\infty}$ $\leq$ $\norm{x}{\infty}$.
Let now $z\in\Li{\infty}{2}$ and $z=x+iy$, with
$x,y\in\Lih{\infty}{2}$; then $\norm{Mz}{\infty}$ $\leq$
$\norm{Mx}{\infty}$ $+$ $\norm{My}{\infty}$ $\leq$
$\norm{x}{\infty}$ $+$ $\norm{y}{\infty}$ $\leq$ $2\norm{z}{\infty}$.
In the same way $\norm{M^*z}{\infty}$ $\leq$
$2\norm{z}{\infty}$, $z\in\Li{\infty}{2}$. \np
 Suppose now $x\in\Li{1}{2}$, as
$B_0 := \{y\in\Li{2}{\infty} : \norm{y}{\infty}\leq1 \}$ is weak$^*$ dense in
the unit ball of $L^\infty$, we have
$$\eqalign{
\norm{Mx}{1} &= \sup\{ |(y,Mx)| : y\in B_0 \} \cr
&= \sup\{ |(M^*y,x)| : y\in B_0 \} \cr
&\leq \sup\{ \norm{M^*y}{\infty} : y\in B_0 \}\norm{x}{1}\leq
2\norm{x}{1}.
}$$
The same holds for $M^*$. \np
 By Riesz-Thorin-Kunze interpolation (proposition 1.8) we have
$\norm{Mx}{p}\leq2\norm{x}{p}$, for $x\in\Li{p}{2}$, and
analogously for $M^*$. \np
 $(ii)$ Let $p\in(1,\infty)$, $x\in L^p_h$, $y\in\Li{1}{\infty}$, and
$e_n:=e_{|x|}(({1\over n},n))$; then $xe_n\to x$ in $L^p$, by
proposition 1.7$(ii)$, and $xe_n\in\Lih{p}{2}$, so that we have \par
$$\eqalign{
(x,\tr{ (M^{(p)}) }y)&=(M^{(p)}x,y)=\lim_{n\to\infty} (M^{(p)}(xe_n),y)
=\lim_{n\to\infty} (M(xe_n),y)\cr
&=\lim_{n\to\infty} (xe_n,M^*y)
=\lim_{n\to\infty} (xe_n,M^{*(p')}y) = (x,M^{*(p')}y).
}$$
Hence the thesis, by linearity and the density of $\Li{1}{\infty}$
in $L^{p'}$. \np
 $(iii)$ Suppose that $p<q<\infty$ and
$x\in\Lih{p}{q}$. Then we have
$$
M^{(p)}x = \lim_{n\to\infty} M^{(p)}(xe_n)
= \lim_{n\to\infty} M^{(q)}(xe_n) = M^{(q)}x,
$$
and, by linearity, we are through. \np
Now consider the case $q=\infty$ and let $x\in\Lih{p}{\infty}$ and $y\in
\Li{1}{\infty}$, then we have
$$\eqalign{
(y,M^{(p)}x)&=\lim_{n\to\infty} (y,M^{(p)}(xe_n))
=\lim_{n\to\infty} (y,M(xe_n)) \cr
&=\lim_{n\to\infty} (M^*y,xe_n)
=(M^{*(1)}y,x)=(y,M^{(\infty)}x).
}$$
So, from proposition 1.7$(iii)$, $M^{(p)}x=M^{(\infty)}x$, and, by
linearity, we are through.
\endproof

\claim {3.4 Corollary} Let $M$ and $M^*$ be sub-Markovian
operators on $L^2$ and $M^{(p)}$, $M^{*(p)}$ their
extensions to $L^p$, $p\in[1,\infty]$. Then $M^{(p)}$, $M^{*(p)}$
are sub-Markovian operators on $L^p$.
\par
\proof{}
Let us consider first the case $p<\infty$. Let $x\in L^p$, $0\leq
x\leq1$, and observe that $x_n:=xe_x(({1\over n},1))\in\Li{2}{p}$, $0\leq
x_n\leq1$ and $x_n\to x$ in $L^p$. Therefore $0\leq M^{(p)}x_n\leq1$ and
$M^{(p)}x_n\to M^{(p)}x$ in $L^p$. From the following lemma it follows
that $0\leq M^{(p)}x\leq1$, that is $M^{(p)}$ is a sub-Markovian operator on
$L^p$. \np
Now consider the case $p=\infty$, and observe that $M^{(\infty)}$ is
obviously weak$^*$ continuous. Besides the ${}^*$-algebra $\Li{2}{\infty}$
is strongly dense in $L^\infty$ so that $\forall x\in L^\infty$, $0\leq
x\leq1$ there exists, by Kaplansky's density theorem, a net
$\{x_\a\}\subset\Li{2}{\infty}$ s.t. $0\leq x_\a\leq1$ and $x_\a\to
x$ strongly hence $\s$-weakly. Then from $M^{(\infty)}x_\a\to
M^{(\infty)}x$ $\s$-weakly and $0\leq M^{(\infty)}x_\a\leq1$ it
follows $0\leq M^{(\infty)}x\leq1$, that is $M^{(\infty)}$ is
sub-Markovian. \np
 A similar proof works also for $M^{*(p)}$.
\endproof

\claim{3.5 Lemma} Let $x\in\L p$, $p\in[1,\infty)$ and $\{x_n\}$ be s.t.
$0\leq x_n\leq1$ and $x_n\to x$ in $L^p$. Then $0\leq x\leq1$
\par
\proof{}
Let $\{x_{n_k}\}$ be s.t. $x_{n_k}\to y\in\L\infty$ weak$^*$, so that
$0\leq y\leq1$. Then $\forall z\in\Li{1}{\infty}$ we get
$(z,y)=\lim_{k\to\infty} (z,x_{n_k})=(z,x)$ and by proposition 1.7$(iii)$
we are through.
\endproof

\claim {3.6 Theorem} Let $\{T_t\}_{t\geq0}$, $\{T^*_t\}_{t\geq0}$
be sub-Markovian semigroups on $L^2$. Then their extensions to $L^p$
are sub-Markovian semigroups on $L^p$, for $p\in[1,\infty]$.
\par
\proof{}
By proposition 3.3, $\norm{T_t}{p}\leq2,\ t\geq0,\ p\in[1,\infty]$. \np
Let $p\in(1,\infty)$, $p'$ the conjugate exponent, $x\in\Li{{p'}}{2}$,
$y\in\Li{p}{2}$, then we have
$$(x,(T_t^{(p)}-I)y)=(x,(T_t-I)y)=((T_t^*-I)x,y)\to0,\ t\to0$$ so that,
from the density of $\Li{{p'}}{2}$ in $L^{p'}$ and of $\Li{p}{2}$ in
$L^p$, we get the weak continuity of $\{T_t\}_{t\geq0}$ on $L^p$, and,
from [D1, proposition 1.23], the strong continuity.
\np
Let now $p=1$ and $e\in Proj(\ca)\cap L^1$, then we have
$$
\norm{T_t^{(1)}e-e}{1}=\norm{T_te-e}{1}\leq
\norm{e}{2}\norm{T_te-e}{2}\to0,\ t\to0.
$$
As $\{x=\sum_{i=1}^n \l_ie_i: \l_i>0,\ e_i\in Proj(\ca)\cap L^1\}$
is total in $L^1$, by proposition 1.7$(i)$, we are through. \np
Now let us observe that the same proof also works for
$\{T^*_t\}_{t\geq0}$, so $\{T^{*(p)}_t\}_{t\geq0}$, $p\in[1,\infty)$, is
a sub-Markovian semigroup on $L^p$. \np
 Finally let $p=\infty$, $x\in L^\infty$ and $y\in\Li{1}{2}$; then we have
$(y,(T^{(\infty)}_t-I)x)=((T_t^{*(1)}-I)y,x)\to0,\ t\to0$, as we have already
proved, so that $\{T^{(\infty)}_t\}_{t\geq0}$ is a weak$^*$ continuous
semigroup
on $L^\infty$. Analogously for $\{T^{*(\infty)}_t\}_{t\geq0}$.
\endproof

Let us now show there is a converse of the preceding theorem in case
$p=\infty$.

\claim {3.7 Theorem} Let $\{T_t\}_{t\geq0}$,
$\{\hat T_t\}_{t\geq0}$ be sub-Markovian semigroups on $\{\ca,\t\}$ s.t.
$$\t(x(T_ty))=\t((\hat T_tx)y),\ x,y\in \Li{1}{\infty}.$$
Then $\{T_t\}_{t\geq0}$ and $\{\hat T_t\}_{t\geq0}$ are the unique
weak$^*$-continuous extensions of sub-Mar\-ko\-vian semigroups on
$L^2$ which are adjoint to each other.  \par
\proof{}
Fix $t\geq0$ and write $T$ for $T_t$ and $\hat T$ for $\hat T_t$.
As $T$ and $\hat T$ are sub-Markovian we get $\norm{T}{\infty}\leq2$,
$\norm{\Th}{\infty}\leq2$. Define $T_*:L^1\to L^1$ and $\Th_*:L^1\to
L^1$ by $\t((\Th_*x)y):=\t(x(Ty))$ and $\t((T_*x)y):=\t(x(\Th y))$,
$x\in L^1,\ y\in L^\infty$. Then $T_*,\ \Th_*$ are positivity
preserving as $x\in L^1_+ \imply \t((\Th_*x)y)=\t(x(Ty))\geq0,\
\forall y\in L^\infty_+$, that is $\Th_*x\geq0$ and analogously
$T_*x\geq0$. Moreover $\norm{\Th_*}{1}=\norm{T}{\infty}$ and
$\norm{T_*}{1}=\norm{\Th}{\infty}$. Besides, $\forall
x,y\in\Li{1}{\infty}$, we get $\t((\Th_*x)y)=\t(x(Ty))=\t((\Th
x)y)$ that is $\Th_*x=\Th x$ and analogously $T_*x=Tx$ for all
$x\in\Li{1}{\infty}$. Therefore, by Riesz-Thorin-Kunze
interpolation (proposition 1.8), $T$ and $\Th$ extend uniquely from
$\Li{1}{\infty}$ to $L^2$ with norm no greater than 2, and
$\t(x(Ty))=\t((\hat Tx)y),\ x,y\in L^2$. \par
Therefore $\{T_t\}_{t\geq0}$ and $\{\hat T_t\}_{t\geq0}$ extend
uniquely to $L^2$. Let now $x,y\in L^2$, then $\forall\e>0$, $\exists\
x',y'\in\Li{1}{\infty}$ s.t. $\norm{x-x'}{2}<\e$, $\norm{y-y'}{2}<\e$
and, as $|(y',(T_t-I)x')|<\e$, for $0\leq t<\d_\e$, we get
$$\eqalign{
|(y,(T_t-I)x)| &\leq |(y-y',(T_t-I)x)|+ |(y',(T_t-I)(x-x'))|+
|(y',(T_t-I)x')| \cr
&\leq \norm{y-y'}{2}\norm{(T_t-I)x}{2}+ |((\Th_t-I)y',x-x')|+
|(y',(T_t-I)x')| \cr
&\leq \norm{y-y'}{2}\norm{(T_t-I)x}{2}+
\norm{(\Th_t-I)y'}{2}\norm{x-x'}{2}+ |(y',(T_t-I)x')| \cr
&\leq \e(3\norm{x}{2}+3(\norm{y}{2}+\e)+1).
}$$
Hence $\{T_t\}_{t\geq0}$ is weakly continuous on $L^2$ and therefore [D1,
proposition 1.23], strongly continuous. Analogously $\{\hat
T_t\}_{t\geq0}$ is strongly continuous on $L^2$. They are
sub-Mar\-ko\-vian semigroups, and the induced extensions of
$T_t|_{\Li{2}{\infty}}$ and $\Th_t|_{\Li{2}{\infty}}$ to $L^\infty$
are the original semigroups.
\endproof

As we saw in proposition 3.3, if $M,M^*$ are sub-Markovian operators on
$L^2$, their extensions to $L^p$ have norm no greater than 2. If we
require a stronger condition of positivity on $M,M^*$ then their
extensions to $L^p$ are contractive operators, as in the following

\claim {3.8 Theorem} If $M$ and $M^*$ are sub-Markovian operators
on $L^2$ s.t.
$$\eqalign{
(Mx)^*(Mx)&\leq\norm{M}{\infty}M(x^*x) \cr
(M^* x)^*(M^* x)&\leq\norm{M^*}{\infty}M^*(x^*x),\quad x\in\Li{2}{\infty},
}$$
then
$$
\norm{M^{(p)}}{p}\leq1,\ \norm{M^{*(p)}}{p}\leq1,
\ p\in[1,\infty].
$$
\par
\proof{}
{}From [DL1, lemma 3.2] one gets
$\norm{M|_{\Li{2}{\infty}}}{\infty}\leq1$ and
$\norm{M^*|_{\Li{2}{\infty}}}{\infty}\leq1$. Let $x\in
L^\infty$ and $\{x_\a\}\subset\Li{2}{\infty}$ s.t. $x_\a\to x$
weak$^*$ and $\norm{x_\a}{\infty}\leq\norm{x}{\infty}$. Then
$M^{(\infty)}x_\a\to M^{(\infty)}x$ weak$^*$ so that
$\norm{M^{(\infty)}x}{\infty}$ $\leq$ $\liminf\norm{M^{(\infty)}x_\a}{\infty}$
$\leq$ $\liminf\norm{x_\a}{\infty}$ $\leq$ $\norm{x}{\infty}$. In the same way
$\norm{M^{*(\infty)}}{\infty}\leq1$. \par
Let now $x\in\Li{1}{2}$, as $B_0:=\{y\in\Li{2}{\infty} :
\norm{y}{\infty}\leq1 \}$ is weak$^*$ dense in the unit ball of
$L^\infty$, we have
$$\eqalign{
\norm{M^{(1)}x}{1}&=\sup\{|(y,M^{(1)}x)| : y\in B_0 \} \cr
&=\sup\{|(M^* y,x)| : y\in B_0 \} \cr
&\leq\sup\{\norm{M^* y}{\infty} : y\in B_0 \}\norm{x}{1} \cr
&\leq\norm{x}{1}
}$$
and, by density and continuity,
$\norm{M^{(1)}x}{1}\leq\norm{x}{1},\ x\in L^1$. By interpolation
$\norm{Mx}{p}\leq\norm{x}{p},\ x\in \Li{2}{p}$ and, by density and continuity,
$\norm{M^{(p)}x}{p}\leq\norm{x}{p},\ x\in L^p$. An analogous result holds
for $M^{*(p)}$.
\endproof

\remark{3.9} The hypotheses of theorem 3.8 are implied by 2-positivity of
the sub-Markovian operators, [DL1]. \par

\claim {3.10 Definition}  Let $\{\ce,\cd(\ce)\}$ be a sesquilinear form on
$L^2(\ca,\t)$. Then
$$
\ce^{[n]}([a_{ij}],[b_{ij}]):=
\sum_{i,j=1}^n \ce(a_{ij},b_{ij}),\ a_{ij},b_{ij}\in \cd(\ce),
$$ is a sesquilinear form on $L^2(\ca\otimes M_n,\t\otimes tr)\cong
L^2(\ca,\t)\otimes L^2(M_n,tr)$, where $tr$ is the usual trace
on $n$ by $n$ matrices. \np
We say that $\ce$ is n-Dirichlet if $\ce^{[n]}$ is a Dirichlet
form.
\par
 \claim{3.11 Lemma}  Let $\{\ce,\cd(\ce)\}$ be a coercive closed form on
$L^2(\ca,\t)$, and let $\{ G_\a \}_{\a\geq0}$ be the associated resolvent.
\np
 Then $\ce^{[n]}$ is a coercive closed form and $\{ G^{[n]}_\a\}_{\a\geq0}$
 is the associated resolvent.
  \par
\proof{}
Let us observe that
\vskip -0.3cm
$$\eqalign{
\ce^{[n]}_\a([a_{ij}],G^{[n]}_\a[b_{ij}]) &=
\ce^{[n]}_\a([a_{ij}],[G_\a b_{ij}]) = \sum_{i,j=1}^n \ce_\a(a_{ij}, G_\a
b_{ij}) \cr
&= \sum_{i,j=1}^n (a_{ij},b_{ij}) = ([a_{ij}],[b_{ij}])
}$$
\vskip -0.3cm\noindent
so that $\{ G^{[n]}_\a\}_{\a\geq0}$ is the resolvent associated to
$\ce^{[n]}$.
Let us now prove that $\{ G^{[n]}_\a\}_{\a\geq0}$ is contractive
\vskip -0.6cm
$$
\norm{\a G^{[n]}_\a[a_{ij}]}{2}^2 = \norm{[\a G_\a
a_{ij}]}{2}^2 \ = \sum_{i,j=1}^n \norm{\a G_\a a_{ij}}{2}^2
\leq \sum_{i,j=1}^n \norm{a_{ij}}{2}^2 = \norm{[a_{ij}]}{2}^2.
$$
\vskip -0.4cm\noindent
Finally let us prove $\{ G^{[n]}_\a\}_{\a\geq0}$ satisfies the sector
condition [see 1.3]. Let $[a_{ij}]$, $[b_{ij}]$  $\in L^2(\ca\otimes
M_n,\t\otimes tr)_h$, then
$$\eqalign{
|([a_{ij}],G^{[n]}_1[b_{ij}])| &= |([a_{ij}],[G_1b_{ij}])| \cr
&= |\sum_{i,j=1}^n (a_{ij},G_1b_{ij})| \cr
&\leq \sum_{i,j=1}^n |(a_{ij},G_1b_{ij})| \cr
&\leq K\sum_{i,j=1}^n (a_{ij},G_1a_{ij})^{1/2}
(b_{ij},G_1b_{ij})^{1/2} \cr
&\leq K(\sum_{i,j=1}^n (a_{ij},G_1a_{ij}))^{1/2}
(\sum_{i,j=1}^n (b_{ij},G_1b_{ij}))^{1/2} \cr
&= K ([a_{ij}],G^{[n]}_1[a_{ij}])^{1/2}
([b_{ij}],G^{[n]}_1[b_{ij}])^{1/2}
}$$
that is $\ce^{[n]}$ is a coercive closed form.
\endproof

\remark{3.12} Analogous results hold for $\{ G^*_\a \}_{\a\geq0}$ and the
associated semigroups $\{T_t\}_{t>0}$ and $\{T^*_t\}_{t>0}$, that is $\{
G^{*[n]}_\a \}_{\a\geq0}$ is the resolvent and $\{T^{[n]}_t\}_{t>0}$, and
$\{T^{*[n]}_t\}_{t>0}$, are the semigroups associated to $\ce^{[n]}$.
\par

\claim {3.13 Theorem}  Let $\{\ce,\cd(\ce)\}$ be a Dirichlet form and $\{
T_t \}_{t\geq0}$, $\{ T^*_t \}_{t\geq0}$ the associated
semigroups. \np
Then $\ce$ is n-Dirichlet $\iff$ $\{T_t\}$ and $\{T^*_t\}$ are
n-positive.
\par
\proof{}
 From the previous lemma it suffices to show that $\{T_t\}$ is sub-Markovian
and $n$-positive $\iff$ $\{T^{[n]}_t\}$ is sub-Markovian and contractive.
  \np
$(\coimply)$ Let $x\in L^2$ be such that $0\leq x\leq1$, then $0\leq x\otimes
1\leq 1\otimes1$, which implies $0\leq T^{[n]}_t(x\otimes1)\leq 1\otimes1$ that
is $0\leq (T_tx)\otimes 1\leq 1\otimes1$ that is $0\leq T_tx\leq1$.  \np
$(\imply)$ First of all, let us observe that $T_t^{(\infty)[n]}$ is the
extension of $T_t^{[n]}$ to $L^\infty(\ca\otimes M_n,\t\otimes tr)$ by
uniqueness, so that $T_t^{(\infty)[n]}$ is positive. \np
Let $x=[x_{ij}]\in L^2(\ca\otimes M_n,\t\otimes tr)$ be s.t. $0\leq x\leq 1$.
Then $0\leq T^{(\infty)[n]}_tx\leq T^{(\infty)[n]}_t1$ that is $0\leq
[T_tx_{ij}]\leq T_t^{(\infty)}1 \otimes 1\leq 1\otimes1$ and the thesis
follows.
\np
Finally $T^{[n]}_t=T_t\otimes 1$ is
obviously a contraction on $L^2(\ca\otimes M_n,\t\otimes tr)$
\endproof

 \section Section 4. Derivations on square integrable operators.
 \par

 In this Section we consider derivations on the space $\L2$. \np
 By this we mean a linear operator
 $$\d:\D\subseteq\L2\to\L2,$$
 where $\D$ is a subalgebra of $\L2\cap\L\infty$, and $\d$ verifies
 $$\d(ab)=a\cdot\d b+\d a\cdot b\qquad a,b\in\D.$$ \par
 We say that a derivation $\d$ is closed under the
${\cal C}^1$, resp. Lipschitz functional calculus if, whenever
$a\in\D_h$, $f(a)\in\D$ for each ${\cal C}^1$, resp. Lipschitz
function $f$ such that $f(0)=0$.
 \par
The domain $\D$ of a derivation is said {\it self-adjoint} if it
is closed under the $^*$ operation. A dense $^*$-subalgebra of
 $\L2\cap\L\infty$ is called a {\it Hilbert algebra}.
 A derivation $\d$ is a {\it $^*$-derivation} if $\D$
is self-adjoint and $\d(a^*)=(\d a)^*$.
 \par
 Now we follow an argument in [Sa] which gives rise to a
{\it non-abelian chain rule} (formula 4.3) for the derivation of
the functional calculus of a self-adjoint element. Let us fix a
self-adjoint element $a\in\A$ and consider the representation
$\pi_a$ of $\C{}\otimes\C{}\equiv{\cal C}_0(\Re\times\Re)$ on $\L2$
given by
 $$
\pi_a(f\otimes g)b=f(a) b g(a)\ ,\qquad b\in\L2.
 $$
and observe that
 $$ Range(\pi_a)\subset\A\vee\A'\ .\eqno(4.1)$$
For each $f\in\Lip$, we set
 $$
\tilde f(s,t)=\cases{
{f(s)-f(t)\over s-t} & $s\not=t$\cr
\vbox{\vskip 2.mm}\cr
f'(t) & $s=t$ }. \eqno(4.2)
 $$
We observe that if $f\in\Lip$ then $\tilde f\in
L^\infty(\Re\times\Re)$, and, if $f(0)=0$, then
$$
\|f\|_\Lip:=\|\tilde f\|_\infty \equiv \norm{f'}{\infty}
$$
is a Banach norm on $\Lip$.
 Now we may state the main theorem of this section:
 \claim{4.1 Theorem} Let $\d$ be a closed derivation on $\L2$, $a\in\cd_h$.
Then the following properties hold:
 \item{$(i)$} $\d$ is closed under the Lipschitz functional calculus.
 \item{$(ii)$} For each $f\in\C1$, $f(0)=0$, one has
 $$\d f(a)=\pi_a(\tilde f)\d a\ .\eqno(4.3)$$
  \item{$(iii)$} For each $f\in\Lip$, $f(0)=0$, one has
$$\|\d f(a)\|_2\leq\|f\|_\Lip\|\d a\|_2\ .\eqno(4.4)$$
 \par
\claim{4.2 Lemma} Let $f$ be a Lipschitz continuous function
such that $f(0)=0$, $\f$ a positive  ${\cal C}^\infty$
function with support in $[-1,1]$
 s.t. $\int\f=1$. Then, the sequence of mollified functions
$\{f_n\}$, $f_n(t):= f*\f_n(t)-f*\f_n(0)$, where
$\f_n(t):=n\f(nt)$, verifies the following properties:
 \item{$(a)$} $f_n(0)=0$
 \item{$(b)$} $\|f-f_n\|_\infty\leq{2\over n}\|f\|_\Lip$
 \item{$(c)$} $\|f_n\|_\Lip\leq\|f\|_\Lip$.
 \item{$(d)$} $\tilde{f_n}\to\tilde f$ weak$^*$ in
$L^\infty(\Re)$
 \par
The proof is trivial and is omitted.
\par
 \claim{4.3 Lemma} Let $\{A,\cd(A)\}$ be a closed linear
operator on $\ch$, $\{x_n\}\subset\cd(A)$ such that
$\|x_n-x\| \to 0$ and exists $k>0$ s.t. $\|Ax_n\| \leq k$. Then there exists
$\{w_n\}$ in the convex hull of $\{x_n\}$ s.t. $w_n\to x$ in the graph-norm of
$A$. As a consequence $x\in\cd(A)$ and $\|Ax\|\leq k$. \par
\proof{}
 Let us denote by $\norm{\cdot}{A}$ the graph-norm,
$\norm{y}{A}:=\|y\|+\|Ay\|$, $y\in\D(A)$. By the hypotheses, the
sequence $x_n$ is bounded in the graph norm, therefore
the set of limit points in the weak topology of
$\{\cd(A),\ \norm{\cdot}{A}\}$ is not empty. By the
Banach-Saks theorem (cf. [DS], Theorem. V.3.14), for any such limit
point $y$, there exists a sequence $w_n$ in the convex hull of $\{x_n\}$  such
that $w_n\to y $ in the graph-norm, hence in the Hilbert norm, therefore $y=x$,
that is $x\in\cd(A)$. Besides
 $$
\|Ax\| = \lim \|Aw_n\| \leq k
 $$
and the thesis follows.
 \endproof
\proof{ of Theorem 4.1} First we observe that, since $a\in\A$,
we may replace $\C1$ by ${\cal C}^1(I)$,
$I:=[-\|a\|,\|a\|]$. Then, equation (4.3) makes sense
also for polynomials, and we check it for $f(t):= t^n$:
 $$\eqalign{\d(a^n)&=\sum_{j=0}^{n-1}a^j(\d a)a^{n-j-1}=\cr
&=\sum_{j=0}^{n-1}\pi_a(s^j t^{n-j-1})\d a\cr
&=\pi_a\left({s^n-t^n\over s-t}\right)\d a
=\pi_a(\tilde f)\d a.\cr}
 $$
By linearity, (4.3) holds for all polynomials $p$ such that
$p(0)=0$.
Finally we observe that, for all such polynomials,
 $$
\eqalign{\|\pi_a(\tilde p)\d a\|_2
&\leq\|\pi_a(\tilde p)\|\|\d a\|_2\cr
&=\| p\|_{Lip(I)}\|\d a\|_2\cr
&=\norm{p}{{\cc^1(I)}} \norm{\d a}{2}.\cr
 }$$
Therefore, if $p_n$ is a sequence of polynomials converging
to a ${\cal C}^1$ function $f$ in the
${\cal C}^1(I)$ norm, then $\d(p_n(a))$ is a
Cauchy sequence w.r.t. the graph norm of $\d$, and $(ii)$
follows by continuity. \np
In particular, we proved that $\d$ is closed under
the ${\cal C}^1$ functional calculus. Therefore, if $f$, $\f_n$ are
as in lemma 4.2, formula (4.3) applies to $f_n$, and we get, using
 lemma $4.2c$,
 $$
\|\d f_n(a)\|_2=\|\pi_a(\tilde f_n)\d a\|_2\leq\|f\|_\Lip\|\d
a\|_2.
 $$
Now we prove that
 $$
\|f(a)-f_n(a)\|_2\to 0.\eqno(4.5)
 $$
Indeed, choosing $\|f-f_n\|_\infty\leq\e^{3}$, we have
 $$
\eqalign{
\|f(a)-f_n(a)\|^2_2
&=\tau \left(\int|f-f_n|^2(\l)de(\l)\right)\leq\cr
&\leq\int_{-\e}^{\e}\left(\left|{f(\l)\over\l}\right|
+\left|{f_n(\l)\over\l}\right|\right)^2d\m(\l)
+{1\over\e^2}\int_{|\l|\geq\e}|f-f_n|^2d\m\leq\cr
&\leq 4\|f\|_\Lip^2\m([-\e,\e])+\|a\|^2_2\e,\cr}
 $$
where the
measure $\m$
is defined by $\m(\Omega):=\tau(\int_\Omega\l^2de(\l))$.
Since $\m(\{0\})=0$ and $\m(\Re)=\|a\|_2$ by definition,
we get $\m([-\e,\e])\to0$, which proves (4.5).
 \np
 Finally we apply lemma 4.3 to the sequence $f_n(a)$
in the domain of $\d$, and $(i)$ and $(iii)$ are proven.
 \endproof
 The following corollary is a consequence of theorem 4.1.
 \claim{4.4 Corollary} Let $\d_n$, $n=1,\dots, N$, be closed
derivations. If $\d:=\sum_n \d_n$, $\D(\d)=\cap_n \D(\d_n)$, then $\d$ is
closed under the Lipschitz functional calculus and, for each $f\in\Lip$,
$f(0)=0$, $$\|\d f(a)\|_2\leq\|f\|_{\Lip}\|\d a\|_2
\qquad\forall a\in\D(\d)_h$$
 \par
\proof{}
By hypothesis $\d$ is closed under Lipschitz functional calculus and (4.3)
holds for $\cc^1$ functions. Let $f$ and $f_n$ be as in lemma 4.2. Then, as
in the proof of theorem 4.1, we get
$$
\norm{\d f_n(a)}{2}=\norm{\pi_a(\tilde f_n)\d a}{2}
\leq \norm{f}{\Lip}\norm{\d a}{2}
$$
and $\norm{f_n(a)-f(a)}{2}\to 0$. Then, by lemma 4.3, we get a sequence
$\{h_n\}$ in the convex hull of $\{f_n\}$  such that $h_n(a)\to f(a)$ in the
graph-norm of $\d_1$. Since $\{h_n(a)\}$ is bounded in the graph-norm of
$\d_2$,
applying again lemma 4.3, we find a sequence in the convex hull
of $\{f_n\}$ converging to $f(a)$ in the graph-norms of $\d_1$ and $\d_2$.
Iterating this procedure $N$ times, we find a sequence $\{g_n\}$ in the convex
hull of $\{f_n\}$ s.t. $g_n(a)$ converges to $f(a)$ in the graph-norms of
$\d_i$, for all $i$. Therefore $g_n(a)\to f(a)$ in the graph-norm of $\d$.
Finally
$$
\norm{\d f(a)}{2}=\lim\norm{\d g_n(a)}{2}
\leq \norm{f}{\Lip}\norm{\d a}{2}.
$$
 \endproof
\rmclaim{4.5 Remark} We would like to compare  theorem 4.1 with an
analogous result of Powers ([Po], cf. theorem 1.6.2 in [B]) for
the derivations on a $C^*$-algebra.
 While the difference in the formulation of the non-abelian
chain-rule (equation 4.3) is just a matter of taste, the difference
on the allowed functional calculus depends on the different norms.
 Indeed, let $\cam$ be the C$^*$-algebra generated by a self-adjoint element in
$\A$. The representation of the tensor product $\cam\otimes\cam$ given by the
left and right actions $\cam$ on $\L2$ extends to a representation of the
$C^*$-tensor product.
 This guarantees that equation (4.3) holds for the closure of the
polynomials in the appropriate norm, $i.e.$ for ${\cal C}^1$
functions.
 On the contrary, if the abelian C$^*$-algebra $\cam$ acts on $\A$, the
tensor product is embedded in the Banach algebra $\cb(\A)$, and
this embedding is not necessarily continuous in the C$^*$ norm. As
a consequence, formula (4.3) does not necessarily hold for ${\cal
C}^1$ functions, as it is shown by McIntosh [Mc].
 \par
 Theorem 4.1 gives a general answer to the problem of the
Lipschitz functional calculus of a self-adjoint operator in
the domain of a derivation. On the other hand there is a
trivial form of ``Lipschitz'' functional calculus that makes
sense also for non self-adjoint operators, $i.e.$ the modulus
of an operator. Now we are going to study this question.
 \par
 \claim{4.6 Lemma} Let $\d$ be a derivation on
$\D$. Then the operator $\d^\dagger:\D^*\to\L2$ defined by
 $$
\d^\dagger a:=(\d a^*)^*\ ,\qquad a\in\D^*
$$ is a derivation. If $\d$ is closed (closable), also
$\d^\dagger$ is.
 \par
 \proof{} The Leibnitz rule for $\d^\dagger$ follows by a
straightforward calculation. The equivalence of the
closability properties follows by the equality
 $$
\|a\|_{\d^\dagger}=\|a^*\|_\d\ ,\qquad\forall a\in\D^*
 $$
 \endproof
 \claim{4.7 Lemma} Let $\d$ be a closed derivation on a
self-adjoint domain $\D$. Then
 $$
\parmat{\d^\dagger&0\cr0&\d\cr}:
\D\otimes M_2\to \L2\otimes M_2
 $$
 is a closed derivation.
 \par
 \proof{} Follows immediately by the property
 $$
\left\|\parmat{\d^\dagger&0\cr0&\d\cr}
\parmat{a&b\cr c&d\cr}\right\|^2=
\|\d a^*\|^2+\|\d b^*\|^2+\|\d c\|^2+\|\d d\|^2.
 $$
 \endproof
 \claim{4.8 Theorem} Let $\d$ be a closed derivation on a
self-adjoint domain $\D$. Then, if $a\in\D$, $|a|\in\D$ and
 $$
\|\d|a|\|_2\leq\sqrt2\|\d a\|_2.\eqno(4.6)
 $$
 \par
\proof{} Consider the functions
 $$\eqalign{
\f_n(t)=
&\sqrt{t+{1\over n^2}}-{1\over n}\qquad t\geq0\cr
\psi_n(t)=
&\sqrt{t^2+{1\over n^2}}-{1\over n}\quad t\in\Re\ .\cr
 }$$
 Since $a\in\D$, $a^*\in\D$, and therefore $a^*a\in\D$ and
 $$\psi_n(|a|)=\f_n(a^*a)\in\D$$
by the theorem on the Lipschitz functional calculus.
\par
Now consider the following chain of inequalities:
 $$
\eqalign{\|\d\psi_n(|a|)\|^2_2
&\leq\|\d\psi_n(|a|)\|^2_2+\|\d\psi_n(|a^*|)\|^2_2=\cr
&=\|\d^\dagger\psi_n(|a|)\|^2_2+\|\d\psi_n(|a^*|)\|^2_2=\cr
&=\left\|\parmat{\d^\dagger&0\cr0&\d\cr}
\parmat{\psi_n(|a|)&0\cr0&\psi_n(|a^*|)\cr}\right\|^2_2=\cr
&=\left\|\parmat{\d^\dagger&0\cr0&\d\cr}
\psi_n\left(\parmat{0&a^*\cr a&0\cr}\right)\right\|^2_2\leq\cr
&=\left\|\parmat{\d^\dagger&0\cr0&\d\cr}
\parmat{0&a^*\cr a&0\cr}\right\|^2_2=\cr
&=\|\d a\|^2_2+\|\d^\dagger a^*\|^2_2=\cr
&=2\|\d a\|^2_2\cr}
 $$
where the main inequality follows by the theorem on
Lipschitz functional calculus applied to the derivation
$\parmat{\d^\dagger&0\cr0&\d\cr}$ mentioned in lemma 4.7.
 \par
 Now it is easy to see that $\psi_n(|a|)\to|a|$ in the
$L^2$ norm, and therefore the result follows from lemma
4.3.
\endproof
We remark that, according to the terminology in [DL1],
theorems 4.1 and 4.8 may be rephrased as follows: each closed
derivation on a Hilbert algebra is a Dirichlet derivation.
 \bigskip
 A natural question related to theorem 4.1 is the
following: when the non-abelian chain rule given in (4.3) extends
to the Lipschitz functional calculus? The first problem is that
$\tilde f$ is not necessarily in the domain of $\pi$, when $f$ is a
Lipschitz function. Indeed $\pi$ may easily be extended to the
C$^*$ tensor product of $L^\infty(\Re)$ with itself, but this space
is smaller than $L^\infty(\Re\times\Re)$ and therefore does not
contain $\tilde f$ in general. Even though we do not try to give
a general answer to the previous question, in the following
proposition we mention two extremal cases in which the addressed
question has a positive answer, the abelian case and the type
I factor case.
 \claim{4.9 Proposition} Let $\A$ be either an abelian algebra or a
type I factor, and $\d$ a closed derivation on $\L2$.
Then, for each self-adjoint $a$ in the domain of $\d$, the map
$\pi_a$ extends to $L^\infty(\Re\times\Re)$. Therefore, the non
abelian chain-rule given in (4.3) extends to Lipschitz functional
calculus.
 \par
 \proof{} It is well known that, if $\A$ is a type I factor, the
map $\pi_a$ extends to a normal representation of
$L^\infty(\s(a)\times\s(a))$, for each self-adjoint $a$ in $\A$.
Then, let $f$ and $f_n$ be as in lemma 4.2. By lemma $4.2d$,
$\pi_a(\tilde{f_n})\d a\to\pi_a(\tilde{f})\d a$ weakly in $\L2$,
and also $\d f_{n_k}(a)\to\d f(a)$ weakly in $\L2$ for a suitable
subsequence, as it is shown in lemma 4.3. Then the
thesis holds.
 If $\A$ is abelian, formula $4.3$ becomes the usual chain rule,
and by normality of the map $f\to f(a)$, $a\in\A_h$, we get the
thesis.
 \endproof
 \rmclaim{4.10 Remark}  We remark that the key property we used in the
proof of the factor case, i.e. the fact that the von~Neumann
algebra generated by the right and left action of any abelian
subalgebra of $\A$ on the standard representation of $\A$ is
isomorphic to the tensor product of $\A$ with itself, is a
characterization of the type I, therefore the property we are
studying is probably confined to type I algebras.
 \par
We conclude this section with an example of a simple
Dirichlet form associated with a general derivation.
\claim{4.11 Proposition} Let $\d$ be a closed derivation
on a Hilbert algebra $\D$. Then the form
 $$\E(x,y)=Re(\d x,\d y)+ Im(\d x,\d y)\eqno(4.7)$$
 is a Dirichlet form.
 \par
 \proof{} The form $\E$ is closed $iff$
its symmetric part $Re(\d x,\d y)$ is, that is $iff$
$\d$ is closed. The weak sector condition follows by
 $$
|\E(x,y)|^2\leq2|(\d x,\d y)|^2\leq2\E(x,x)\E(y,y).
 $$
 Now let us define the operators
 $$
\eqalign{d_1a=&{\d a+\d^\dagger a\over 2}\cr
         d_2a=&{\d a-\d^\dagger a\over 2i}\cr}\qquad\qquad a\in\D.
 $$
 It is clear that $d_1$, $d_2$ are $^*$-derivations and
 $$
\E(x,y)=\sum_{i,j\in\{1,2\}}a_{ij}(d_ix,d_jy),
 $$
where $(a_{ij})=\left(\matrix{1&1\cr-1&1\cr}\right)$.
Now Dirichlet property follows from the general theorem in
the following section.
\endproof
As the proof of the preceding proposition shows, the study of the
Dirichlet form in (4.7) may be reduced to the case of
$^*$-derivations. Therefore in the following section only
$^*$-derivations will be considered.

\section Section 5. Explicit constructions of Dirichlet forms.
\par

The aim of this section is to describe a class of Dirichlet forms which can be
considered as the non commutative generalization of a class of
(generally non symmetric) commutative Dirichlet forms studied in [MR].
At the same time these non commutative examples also extend previous ones
constructed in [DL1]. \par
We start considering the following leading example taken from the commutative
context. Let $A=[a_{ij}]$ be an element of $L^1_{loc}(U)\otimes M_n$,
$U\subset{\Re^n}$ open. Then we define the bilinear form on $\cc^\infty_o(U)$
 $$
\E(u,v):=\sum_{i,j=1}^{n}\int
a_{ij}(x){\partial{u}\over\partial{x_i}}
{\partial{v}\over\partial{x_j}}dx\ .\eqno(5.1)
$$
Under the following simple assumptions: there exists $0<\nu <\infty$ such that
$$\cases{
\sum_{i,j=1}^{n}\tilde a_{ij}(x)\xi_{i}\xi_{j}\geq\nu{\|\xi\|}^2\quad
&\qquad $\forall\xi=(\xi_{1},\dots\xi_{n})$\cr
\vbox{\vskip2.0mm}\cr
\check a_{ij}\in L^\infty(U,dx)\quad&\qquad $ 1\leq i,j\leq n $.\cr}
\eqno(5.2)
$$
where $\tilde{A}$ (resp.$\check A$) is the symmetric (resp. antisymmetric) part
of $A$, it can be proven that the form $\{\E,D(\E)\}$ is closable and its
closure is a Dirichlet form (cf.[MR]).
\par
Now we discuss some generalizations of the
preceding example to the non commutative context. The general setting is the
following: we have a sesquilinear form of the type
 $$
\E(x,y):=\sum_{i=1}^n(d_ix,a_{ij}d_jy)\eqno(5.1')
 $$
where $d_i$, $i=1\dots n$, is a family of $^*$-derivations on $\L2$ and the
$a_{ij}$'s belong to the center of $\A$. In this case condition (5.2) is
replaced by its proper non-commutative analogue:
$$\cases{
\tilde A\geq\nu I\ ,&\cr
\vbox{\vskip2.0mm}\cr
\check a_{ij}\in\L\infty\ ,&\quad $ 1\leq i,j\leq n $.\cr}
\eqno(5.2')
$$
 \par
 Our first result is theorem 5.1. This theorem  deals with {\it
general} $^*$-derivations with a dense common domain. In this case, because of
such generality, we need stronger requirements in order to get  closedness for
the form $(5.1')$. This is obtained by asking the symmetric part of the matrix
$A$ to be the identity matrix, which makes the first condition in $(5.2')$
automatically fulfilled.
 \par
 In theorem 5.2 the derivations $d_i$ are given by commutators $[z_i,\cdot]$,
where the $z_i$'s  are skew-symmetric elements in $L^2+L^\infty$, which
provides
the closedness of the form $\E$. In this case conditions $(5.2')$ suffice to
get
the result. Such a theorem is a non symmetric extension of theorem 6.10 in
[DL1].
 \claim{5.1 Theorem} Assume we are given a family
$$
d_{i}:\D_i\subset\L2\to\L2 \quad i=1,...,n
$$
of $^*$-derivations over Hilbert algebras ${\D_i}$ such that
\item{$(a)$} each $d_i$ is closable
\item{$(b)$} $\D:=\cap_{i=1}^{n}\D_i$ is dense
\np
and consider the form  $\E$ given by
 $$
\left\{\eqalign{
D(\E)  &:=\D                                                     \cr
\E(x,y)&:=\sum_{i=1}^n(d_ix,d_iy)+\sum_{i,j=1}^n(d_ix,c_{ij}d_jy)\cr
}\right.\eqno(5.3)
 $$
where the $c_{ij}$'s are self-adjoint elements in the center of $\A$ such that
$c_{ij}=-c_{ji}$. Then the form is closable and its closure is a Dirichlet
form.
\par
\proof{} Sesquilinearity of $\E$ being evident, we prove real positivity.
We have
$$
\E(x,x)=\sum_{i=1}^{n}\|d_{i}x\|^{2}+\t\otimes tr(C B(x,x))\eqno(5.4)
$$
with $C=[c_{ij}]\in \L\infty\otimes M_n$ and
$B(x,y)\in \L1\otimes M_n$ is given by
$B(x,y)=[d_{i}yd_{j}x]$. Since $C$ is a real antisymmetric matrix and
$B(x,x)$ is a real symmetric matrix when $x\in\D_h$, the last term in the right
hand side of (5.4) vanishes, and the positivity of $\E$ follows.
 Because of hypothesis $(a)$ the form $\tilde\E$ is closable (see e.g. [DL1])
and therefore, by definition, $\E$ is closable. We now prove that the weak
sector condition holds for $\E$, $i.e.$ there exists  $0<K<\infty$ such that
$$
|\check{\E}(x,y)|\leq
K{\tilde{\E}_{1}(x,x)}^{1/2}{\tilde{\E}_{1}(y,y)}^{1/2} $$
for all $x,y\in\D_h$, where $\tilde{\E}_{1}(x,y):=\tilde{\E}(x,y)+(x,y)$.\par
Setting $M:=\|C\|_\infty$ and applying H\"older and Schwartz inequalities we
get
$$
\eqalign{
|\check\E(x,y)| &= |\t\otimes tr(C B(x,y))|
\leq M\|B(x,y)\|_1\cr
&\leq M\sum_{i,j=1}^n\|B(x,y)_{ij}\|_1
\leq M\sum_{i,j=1}^{n}\|d_{j}x\|_2\|d_{i}y\|_2 \cr
&\leq nM (\sum_{j=1}^{n}{\|d_{j}x\|}_2^2)^{1/2}
(\sum_{i=1}^{n}({\|d_{i}y\|}_2^2) ^{1/2}
\leq nM {\tilde\E_{1}(x,x)}^{1/2}
{\tilde\E_{1}(y,y)}^{1/2}.\cr
}
$$
It remains to prove that the
closure of the form (5.3) is a Dirichlet form. Let $\overline{\E}$ denote the
closure of the form $\E$, which is obtained replacing the $d_i$'s with their
closures $\ov d_i$ (cf. [DL1]).
 We notice that if $x\in\D(\overline\E)$ then
$\varphi_{0}(x):= x^+\wedge1\in\D(\overline\E)$ because this holds for each
$\D_{i}$. Hence we claim that $$
\overline{\E}(x\mp\varphi_{0}(x),x\pm\varphi_{0}(x))\geq
0\quad\forall x\in\D(\overline\E)_h\ .\eqno(5.5)  $$
 First we observe that the matrix $[(\ov d_i x,c_{ij}\ov d_j\f_0(x))]$ is
antisymmetric. It is enough to show this replacing $\f_0$ with a smooth
approximation $\f$. Then, by equations (4.1) and (4.3) and by the
hypotheses on $C$, we get
 $$
\eqalign{
(\ov d_i x,c_{ij}\ov d_j\f(x))&=(\ov d_i x,c_{ij}\pi(\tilde\f)\ov d_jx)\cr
&=(\pi(\tilde\f)\ov d_i x,c_{ij}\ov d_jx)\cr
&=(\ov d_i \f(x),c_{ij}\ov d_jx)\cr
&=-(\ov d_j x,c_{ji}\ov d_i\f(x)).\cr}
 $$
Therefore, using again the antisymmetry of $C$, we have
 $$
\overline{\E}(x\mp\varphi_{0}(x),x\pm\varphi_{0}(x))
=\sum_{i=1}^{n}(\ov d_{i}(x\mp\varphi_{0}(x)),\ov d_{i}(x\pm\varphi_{0}(x)))
 $$
If we set  $B^\mp_{ij}:= (\ov d_{i}(x\mp\varphi_{0}(x)),\ov
d_j(x\pm\varphi_0(x)))$, the left hand side in (5.5) is just the trace of
$B^{\mp}$, therefore (5.5) follows if we prove that $B^{\mp}$ is a positive
definite real-valued matrix. Indeed, for any $(t_{1},\dots,t_{n})\in\Re^n$ and
setting $d:=\sum_1^n \ov d_i$, we have
 $$
\eqalign{
\sum_{i,j=1}^nB^\mp_{ij}t_{i}t_{j}
&=((\sum_{i=1}^{n}t_{i}\ov d_{i})(x\mp\varphi_{0}(x)),
(\sum_{j=1}^{n}t_{i}\ov d_{j})(x\pm\varphi_{0}(x)))\cr
&=(d(x\mp\varphi_{0}(x)),d(x\pm\varphi_{0}(x)))\cr
&=(dx,dx)-(d\varphi_{0}(x),d\varphi_{0}(x))\geq 0\cr}
 $$
by corollary 4.4, and this ends the proof.
\endproof
\claim{5.2 Theorem}
Let $z_{1},\dots,z_{n}$ be skew-adjoint elements in
$L^2+L^\infty$, define
$$
d_{i}(x):=z_{i}x-xz_{i}\quad\forall x\in \Li2\infty
$$
and let $A=[a_{ij}]$ be a matrix of self-adjoint elements in the
center of $\A$ such that condition $(5.2')$ holds. Then, the form
$$
\left\{\eqalign{
\D(\E) &:= \Li2\infty\cr
\E(x,y)&:= \sum_{i,j=1}^n (d_i x,a_{ij}d_j y)\cr
}\right.
$$
is closable and its closure is a Dirichlet form.
\par
\proof{} Let us denote by $B$ the square root of the symmetric part of $A$
as an element in $\L\infty\otimes M_n$. We also set $\d_i:=\sum b_{ij}d_j$ and
$C:=B^{-1}\check A B^{-1}$.  We notice that, since $\tilde A$ is
coercive, $B^{-1}$ is bounded and $C$ is bounded, real and skew-symmetric.
 Then,
 $$
\E(x,y)=\sum_{i=1}^n(\d_ix,\d_jy)+(\d_ix,c_{ij}\d_jy).
 $$
Since $z_i\in L^2+L^\infty$, $i=1,\ldots,n$ and $B$ is bounded, $w_i:=\sum
b_{ij}z_j\in L^2+L^\infty$. Then $\d_i$, being implemented by $w_i$, is a
closable derivation on $\Li2\infty$ (see e.g. [DL1]) and the thesis follows by
theorem 5.1.
 \endproof

 \section Acknowledgments. \par

The authors would like to thank S.~Albeverio for
suggesting them this line of research, and R.~Longo for discussions. They would
also like to thank the organizers of  the Nottingham Conference on Quantum
Probability (March 93) where the results of this paper were announced.

 \section References. \par

\item{[AH]} S. Albeverio, R. H\o egh-Krohn. {\it Dirichlet forms and Markov
semigroups on $C^*$-algebras.} Comm. Math. Phys. {\bf 56} (1977) 173.

\item{[AHO]} S. Albeverio, R. H\o egh-Krohn, G. Olsen.
{\it Dynamical semigroups and Markov processes
on $C^*$-algebras.} Journ. f\"ur Math. {\bf 319} (1980) 25.

\item{[AMR]} S. Albeverio, Z.M. Ma, M. R\"ockner. {\it Quasi-regular Dirichlet
Forms and Markov Processes.} Journ. Funct. An. {\bf 111} (1993) 118.

\item{[AR]} S. Albeverio, M. R\"ockner. {\it Classical Dirichlet Forms on
topological vector spaces - Construction of an associated diffusion processes.}
Probab. Theory Related Fields {\bf 83} (1989) 405.

\item{[AW]} L. Accardi, W. von Waldenfels. {\it Quantum probability and
applications.} Lecture Notes in Math. 1136 (1985), 1303 (1988), 1396 (1989),
1442 (1990), Springer, New York.

\item{[B]} O. Bratteli. {\it Derivations, Dissipations and Group actions on
$C^*$-algebras.} Lecture Notes in Math. 1229. Springer, New York 1986.

%\item{[BR]} O. Bratteli, D.W. Robinson. {\it Operator algebras
%and quantum   statistical mechanics, I.} Sprin\-ger, New York, 1979.

\item{[Ci]} F. Cipriani. {\it Dirichlet forms and Markovian semigroups on
standard forms of von Neumann algebras.} SISSA PhD Thesis 1992.

\item{[Co]} A. Connes. {\it Geometrie non commutative.} Intereditions, Paris,
1990.

\item{[D1]} E.B. Davies. {\it One-parameter semigroups.} Academic Press, 1980.

\item{[D2]} E.B. Davies. {\it Heat kernels and spectral theory.} Cambridge
University Press, 1989.

\item{[D3]} E.B. Davies. {\it Analysis on Graphs and Non Commutative Geometry.}
Journ. Funct. An. {\bf 111} (1993) 398.

%\item{[Di]} J. Dixmier. {\it Von Neumann Algebras.} North Holland 1981.

\item{[DL1]} E.B. Davies, J.M. Lindsay. {\it Non-commutative symmetric Markov
semi\-groups.} Math. Z. {\bf 210} (1992) 379.

\item{[DL2]} E.B. Davies, J.M. Lindsay. {\it Superderivations and symmetric
Markov semigroups.} To appear on Comm. Math. Phys.

\item{[DR]} E.B. Davies, O.S. Rothaus. {\it Markov semigroups on
$C^*$-bundles}. Journ. Funct. An. {\bf 85} (1989) 264.

\item{[DS]} N. Dunford, J.T. Schwartz. {\it Linear operators. Part I : General
Theory.} Interscience Publ., New York, 1964.

\item{[F]} M. Fukushima. {\it Dirichlet forms and Markov processes}. North
Holland, 1980.

\item{[G]} L. Gross. {\it Hypercontractivity and logarithmic Sobolev
inequalities for the Clifford Dirichlet forms.} Duke Math. J. (3) {\bf 42}
(1975) 383.

\item{[GL]} S. Goldstein, J.M. Lindsay. {\it KMS-symmetric Markov semigroups.}
(preprint) ; {\it Beurling-Deny conditions for KMS-symmetric dynamical
semigroups.} (preprint).

\item{[Mc]} A. Mc Intosh. {\it Functions and derivations on $C^*$-algebras.}
Journ. Funct. An. {\bf 30} (1978) 264.

\item{[MR]} Z.M. Ma, M. R\"ockner. {\it An introduction to the theory of
(non-symmetric) Dirichlet forms.} Lecture Notes in Math. Springer, New
York, 1992.

\item{[N]} E. Nelson. {\it Notes on non-commutative integration.} Journ.
Funct. An. {\bf15} (1974) 103.

\item{[Pa]} K.R. Parthasarathy. {\it An introduction to quantum stochastic
calculus.} Birkh\"auser, 1992.

\item{[Po]} R.T. Powers. {\it A remark on the domain of an unbounded
derivation of a $C^*$-algebra.} Journ. Funct. An. {\bf18} (1975) 85.

%\item{[Sak]} S. Sakai. {\it $C^*$-Algebras and $W^*$-Algebras.} Springer 1971.

\item{[Sa]} J.L. Sauvageot. {\it Quantum Dirichlet forms, differential calculus
and semigroups.} In: Quantum Probability and applications V. Proceedings
Heidelberg 1988, LNM 1442 Springer, New York.

\item{[Se]} I.E. Segal. {\it A non-commutative extension of abstract
integration.} Ann. Math {\bf57} (1953) 401.

\item{[St]} W.F. Stinespring. {\it Integration theorems for gages and duality
for unimodular groups}. Trans. Am. Math. Soc. {\bf 90} (1959) 15.

%\item{[Str]} S. Stratila. {\it Modular theory in operator algebras.} Abacus
%Press, Turnbridge Wells, 1981.

\item{[T]} M. Takesaki. {\it Theory of operator algebras, I.} Springer,
New York, 1979.

\end